\documentclass[12pt]{article}
\usepackage{amsmath,amssymb,amsthm}
\title{Rigid systems of
second-order linear
differential equations\footnotetext{This is the authors' version of a work that was published in Linear Algebra Appl. 414 (2006) 517--532.}}

\author{
M$^{\text{\b a}}$  Isabel
Garc{\'\i}a-Planas\\
 Universitat Polit\`ecnica de Catalunya\\
 Miner{\'\i}a 1, Esc. C, 1-3,
 08038 Barcelona, Spain
\\
maria.isabel.garcia@upc.edu,\and
M. Dolors Magret \\
 Universitat Polit\`ecnica de Catalunya\\
 Miner{\'\i}a 1, Esc. C, 1-3,
 08038 Barcelona, Spain
\\
m.dolors.magret@upc.edu,\and
Vladimir V. Sergeichuk%
\thanks{Corresponding author}\\ Institute of
Mathematics\\
Tereshchenkivska St.
3, Kiev, Ukraine\\
sergeich@imath.kiev.ua\and
Nadya A. Zharko\\
Mech.-Math. Faculty,
Kiev National University,\\
Vladimirskaya 64, Kiev, Ukraine\\
n.zharko@mail.ru}

\date{}

\begin{document}
\newcommand{\rank}{\mathop{\rm rank}\nolimits}
\newcommand{\diag}{\mathop{\rm diag}\nolimits}
\newcommand{\im}{\mathop{\rm Im}\nolimits}

\newtheorem{theorem}{Theorem}
\newtheorem{lemma}[theorem]{Lemma}
\newtheorem{corollary}[theorem]{Corollary}

\theoremstyle{definition}
\newtheorem{definition}[theorem]{Definition}

\theoremstyle{remark}
\newtheorem{example}[theorem]{Example}
\newtheorem{remark}[theorem]{Remark}

\renewcommand{\le}{\leqslant}
\renewcommand{\ge}{\geqslant}

\maketitle

\begin{abstract}
We say that a system of differential
equations
\[
 \ddot{x}(t)=A\dot
 x(t)+Bx(t)+Cu(t),
\qquad A,B\in{\mathbb C}^{m\times
m},\ C\in{\mathbb C}^{m\times n},
\]
is rigid if it can be reduced by
substitutions
\[
 x(t)=Sy(t),\qquad
u(t)=U\dot y(t)+Vy(t)+Pv(t)
\]
with nonsingular $S$ and $P$ to each
system obtained from it by a small
enough perturbation of its matrices
$A,B,C$. We prove that there exists
a rigid system if and only if
$m<n(1+\sqrt{5})/2$, and describe
all rigid systems.


\end{abstract}

\section{Introduction}   \label{s1}

We consider a system of
differential equations
\begin{equation}\label{1.5}
 \ddot{x}(t)=A\dot x(t)+Bx(t)+Cu(t),
\qquad A,B\in{\mathbb
C}^{m\times m},\ C\in{\mathbb
C}^{m\times n},
\end{equation}
in which $x(t)$ is the
unknown vector function,
$u(t)$ is a vector function,
and $\dot x(t)=d x(t)/dt$.
Any substitution
\begin{align*}\label{1.2}
 x(t)=&Sy(t),\\
u(t)=&U\dot y(t)+Vy(t)+Pv(t)
\end{align*}
with nonsingular $S$ and $P$
transforms it to the system
\[
\ddot{y}(t)=S^{-1}(AS+CU)\,\dot
y(t) +S^{-1}(BS+CV)\,y(t)
+S^{-1}CP\,v(t),
\]
which has the form
\eqref{1.5} and is given by
the matrices
\begin{equation*}\label{rtg}
A'=S^{-1}(AS+CU),\qquad
B'=S^{-1}(BS+CV),\qquad
C'=S^{-1}CP.
\end{equation*}
In partitioned matrix
notation
\begin{equation}\label{1.4b}
  [C'\ B'\ A']=
  S^{-1}[C\ B\ A]\begin{bmatrix}
    P&V&U\\0&S&0\\0&0&S
  \end{bmatrix}.
\end{equation}

\begin{definition}
By an $m\times (n,m,m)$
triple we mean a triple of
${m\times n}$, ${m\times m}$,
and ${m\times m}$ matrices.
Two such triples $(C,B,A)$
and $(C',B',A')$ are said to
be \emph{feedback similar} if
they satisfy \eqref{1.4b} for
some $V$, $U$, and
nonsingular $P$ and $S$. (The
term ``feedback similarity''
comes from systems theory.)
\end{definition}

Every transformation of
feedback similarity on a
triple $(C,B,A)$ can be
realized by a sequence of the
following operations:
\begin{itemize}
  \item[(i)]
A simultaneous elementary row
operation on $C,\ B$, and
$A$, and then the inverse
column operation on $B$ and
the inverse column operation
on $A$.

 \item[(ii)]
An elementary column
operation on $C$.

 \item[(iii)]
Adding any constant multiple
of a column of $C$ to a
column of $B$ or $A$.
\end{itemize}

The matrices $A$, $B$, and
$C$ are written in the block
matrix $[C\ B\ A]$ in the
reverse order to ensure that
all admissible additions of
columns are performed from a
left block to a right block
as is customary in matrix
problems (see, for instance,
\cite{gab+roi} or
\cite{ser}).

Related matrix problems are
considered by systems
theorists \cite{glu, hin1,
hin2, loi}.

The canonical form problem
for a matrix triple $(C,B,A)$
up to feedback similarity is
hopeless even if $C=0$ since
then the pair $(B,A)$ reduces
by simultaneous similarity
transformations, and the
problem of classifying pairs
of matrices up to
simultaneous similarity
contains both the problem of
classifying \emph{any} system
of linear operators and the
problem of classifying
representations of \emph{any}
finite-dimensional algebra
\cite{bel-ser}.
Classification problems that
contain the problem of
classifying matrix pairs of
up to simultaneous similarity
are called \emph{wild}.

Nevertheless, using Belitskii's
algorithm \cite{bel,ser} one can
reduce any given triple ${\cal
T}=(C,B,A)$  by transformations
(i)--(iii) to some \emph{canonical}
triple ${\cal T}_{\text{can}}$; this
means that ${\cal T}_{\text{can}}$
is feedback similar to ${\cal T}$
and two triples ${\cal T}$ and
${\cal T}^{\,\prime}$ are reduced by
Belitskii's algorithm to the same
triple ${\cal T}_{\text{can}}={\cal
T}_{\text{can}}^{\,\prime}$ if and
only if ${\cal T}$ and ${\cal
T}^{\,\prime}$ are feedback similar.
(Of course, an explicit description
of \emph{all} canonical matrices
does not exist since the matrix
problem is wild.)

A canonical form problem
simplifies if the matrices
are considered up to
arbitrarily small
perturbations (this case is
important for applications in
which one has matrices that
arise from physical
measurements since then their
entries are known only
approximately). For instance,
a square matrix $A$ reduces
to a diagonal matrix $D$ by
an arbitrarily small
perturbation (making its
eigenvalues pairwise
distinct) and similarity
transformations. The matrix
$D$ is determined by $A$ up
to small perturbations of
diagonal entries.

In Lemma \ref{t4.1} we give a
normal form of $m\times
(n,m,m)$ triples for
arbitrarily small
perturbations and feedback
similarity. A canonical form
of such triples if $n$
divides $m$ is obtained in
Theorem \ref{c4.2}.

By analogy with quiver
representations \cite[p.
203]{der-wey}, we say that a
matrix $t$-tuple $\cal A$ is
\emph{rigid} with respect to
some equivalence relation on
the set of $t$-tuples of the
same size if there is a
neighborhood of $\cal A$
consisting of $t$-tuples that
are equivalent to $\cal A$.
For instance, the matrices
$I,\ [I\; 0]$, and $[I\;
0]^T$ are rigid with respect
to elementary
transformations, but each
matrix is not rigid with
respect to similarity
transformations.

In Theorem \ref{c4.1} we
prove that there exists an
$m\times (n,m,m)$ triple that
is rigid with respect to
feedback similarity if and
only if
\begin{equation}\label{yuf}
m<\frac {1+\sqrt{5}}2\,n.
\end{equation}
We also construct such a
rigid triple ${\cal T}_{mn}$
for each $m$ and $n$
satisfying \eqref{yuf} and
prove that each $m\times
(n,m,m)$ triple reduces to
${\cal T}_{mn}$ by an
arbitrarily small
perturbation and a feedback
similarity transformation (so
${\cal T}_{mn}$ can be
considered as a canonical
triple for arbitrarily small
perturbations and feedback
similarity). All triples that
reduce to ${\cal T}_{mn}$ by
feedback similarity
transformations form an open
and everywhere dense set in
the space of all $m\times
(n,m,m)$ triples; moreover,
this set consists of all
rigid triples of this size.

The mentioned results about
triples will be obtained in
Section \ref{s4}.

In Section \ref{s3} we
consider analogous problems
for systems of first-order
linear differential
equations. Such a system is
given by a matrix pair; the
results of Section \ref{s3}
are used in Section \ref{s4}.

In Section \ref{s1aa} we
prove a technical lemma.

\section{Perturbations}
\label{s1aa}

The \emph{norm} of a complex matrix
$A=[a_{ij}]$ is the nonnegative real
number
\[
\|A\|=\sqrt{\sum|a_{ij}|^2}.
\]
For each $m\times (n,m,m)$
triple ${\cal P}=(C,B,A)$, we
denote
\begin{equation*}\label{erf}
\|{\cal P}\|:=
\|C\|+\|B\|+\|A\|
\end{equation*}
and define the block matrix
\[[{\cal P}]:=[C\ B\ A].
\]

We say that a matrix triple
${\widetilde {\cal T}}$ is
obtained from ${\cal T}$ by a
sequence of perturbations and
feedback similarity
transformations if there is a
sequence of triples
\begin{equation*}\label{fr}
 {\cal T}={\cal T}_1,\ \ {\cal
 T}_2,\ \
{\cal T}_3,\dots,\ {\cal
T}_{l+1}=\widetilde{\cal F},
\end{equation*}
in which
\begin{equation}\label{nmj}
{\cal T}_2=S_1^{-1}[{\cal T}+
 \Delta{\cal
 T}_1]R_1,\quad {\cal T}_3
 =S_2^{-1}[{\cal T}_2+
 \Delta {\cal
 T}_2]R_2,\ \ldots
\end{equation}
($\Delta {\cal T}_1, \Delta
{\cal T}_2,\dots$ are
triples), and all $R_i$ have
the form \eqref{1.4b}:
\begin{equation*}
 R_i=\begin{bmatrix}
    P_i&V_i&U_i\\0&S_i&0\\0&0&S_i
  \end{bmatrix}.
\end{equation*}

\begin{lemma} \label{lemm2}
Let $\varepsilon$ be any
positive number and let a
triple ${\widetilde {\cal
T}}$ be obtained from a
triple ${\cal T}$ by a
sequence \eqref{nmj} of
perturbations and feedback
similarity transformations
satisfying
\begin{equation*}\label{2.6a}
\|\Delta {\cal T}_1\|<
\frac{\varepsilon}{2}\,,\qquad
\|\Delta{\cal T}_{i+1}\|<
\frac{\varepsilon}{2^{i+1}\|\widetilde
S_{i}\| \|\widetilde
R_{i}^{-1}\|} \quad
(i=1,2,\dots,l),
\end{equation*}
where
 \[
 \widetilde S_i:=S_1S_2\cdots S_i,\qquad
\widetilde R_i:=R_1R_2\cdots
R_i.
\]
Then ${\widetilde {\cal T}}$
is feedback similar to some
triple ${\cal T}+\nabla{\cal
T}$, $\|\nabla {\cal T}\|<
\varepsilon$.
\end{lemma}

\begin{proof}
If $l=2$, then by \eqref{nmj}
\begin{align*}
[\widetilde{\cal F}]&=[{\cal
T}_3]
 =S_2^{-1}[{\cal T}_2+
 \Delta {\cal T}_2]R_2=
 S_2^{-1}(S_1^{-1}[{\cal T}+
 \Delta {\cal T}_1]R_1+
 [\Delta{\cal T}_2])R_2\\
 &=
 (S_1S_2)^{-1}([{\cal T}]+
 [{\Delta {\cal T}_1}]+
 S_1[\Delta {\cal T}_2]
 R_1^{-1})R_1R_2.
\end{align*}
Analogously, for any $l$
\[
 [\widetilde{\cal F}]=[{\cal T}_{l+1}]=
 \widetilde S_l^{-1}[{\cal
 T}+\nabla {\cal T}]\widetilde R_l,
\]
where
$$
 [\nabla {\cal T}]:=
 [\Delta{\cal T}_1]+
 \widetilde S_1[\Delta {\cal
 T}_2]
 \widetilde R_1^{-1}+\dots+
 \widetilde S_{l-1}[\Delta
 {\cal T}_l] \widetilde
 R_{l-1}^{-1}.
$$
Then
\begin{align*}
 \|\nabla {\cal T}\|&\le
 \|\Delta{\cal T}_1\|+
 \|\widetilde S_1\|\cdot
 \|\Delta {\cal T}_2\|\cdot
 \|\widetilde R_1^{-1}\|+\dots+
 \|\widetilde S_{l-1}\|\cdot
 \|\Delta
 {\cal T}_l\|\cdot \| \widetilde
 R_{l-1}^{-1}\|
 \\
 &<\frac {\varepsilon}2 +
 \frac {\varepsilon}4 +\dots+
 \frac {\varepsilon}{2^l} < \varepsilon .
\end{align*}
\end{proof}

\begin{corollary}\label{corrd}
Let a matrix triple ${\cal
T}$ reduce to a triple from
some set $\cal S$ by a
sequence of arbitrarily small
perturbations and feedback
similarity transformations.
Then ${\cal T}$ is
transformed by some
arbitrarily small
perturbation to a triple that
is feedback similar to a
triple in $\cal S$.
\end{corollary}

\section{Feedback similarity
of matrix pairs} \label{s3}

In this preliminary section
we consider problems studied
in Section \ref{s4} in much
simpler case: for systems of
first-order linear
differential equations
\begin{equation}\label{a1}
\dot x(t)=Ax(t)+Bu(t), \qquad
A\in{\mathbb C}^{m\times m},\
B\in{\mathbb C}^{m\times
n}.
\end{equation}
Any substitution
\begin{equation*}\label{1.2cf}
\begin{array}{l}
 x(t)=Sy(t),\\
u(t)=Uy(t)+Pv(t)
\end{array}
\end{equation*}
with nonsingular $S$ and $P$
transforms it to the system
\begin{equation*}\label{ji}
\dot{y}(t)=S^{-1}(AS+BU)\,
y(t) +S^{-1}BP\,v(t)
\end{equation*}
of the form \eqref{a1}, whose
matrices $A'$ and $B'$ can be
calculated as follows:
\begin{equation}\label{1.4a}
  [B'\ A']=
  S^{-1}[B\ A]\begin{bmatrix}
    P&U\\0&S
  \end{bmatrix}.
\end{equation}

In systems theory, \eqref{a1}
is called the standard linear
system (without output), $A$
is called the system matrix,
$B$ is called the input
matrix, $u(t)$ is the input
to the system at time $t$ (it
is the way that the external
world affects the system),
and $x(t)$ is the state of
system at time $t$ (it is the
memory of the net effect of
past inputs). The system
\eqref{a1} is said to be
\emph{controllable} if the
spectrum of $A+BU$ can be
placed arbitrarily by choice
of $U$, this holds if and
only if $\rank [B\ AB\
\ldots\ A^mB]= m$.

\begin{definition}
By an $m\times (n,m)$ pair we
mean a pair of $m\times n$
and $m\times m$ matrices. Two
such pairs $(B,A)$ and
$(B',A')$ are said to be
\emph{feedback similar} if
they satisfy \eqref{1.4a} for
some $U$ and nonsingular $P$
and $S$.
\end{definition}

Every feedback similarity
transformation on $(B,A)$ can
be realized by a sequence of
the following operations:
\begin{itemize}
  \item[(i$'$)]
A simultaneous elementary row
operation on both matrices,
and then the inverse column
operation on $A$.

 \item[(ii$'$)]
An elementary column
operation on $B$.

 \item[(iii$'$)]
Adding any constant multiple
of a column of $B$ to a
column of $A$.
\end{itemize}

In the next section we will
reduce a triple $(C,B,A)$ to
canonical form for
arbitrarily small
perturbations and feedback
similarity using results of
this section as follows.
First we reduce its subpair
$(C,B)$ to the pair $(C_{\rm
can},B_{\rm can})$ defined in
\eqref{3.2a}, which is
canonical with respect to
arbitrarily small
perturbations and feedback
similarity; respectively,
$(C,B,A)$ reduces to some
triple $(C_{\rm can},B_{\rm
can}, A')$. Then we reduce
$A'$ to canonical form for
arbitrarily small
perturbations and those
feedback similarity
transformations on $(C_{\rm
can},B_{\rm can}, A')$ that
preserve $(C_{\rm can},B_{\rm
can})$, these transformations
are described in Theorem
\ref{t3.2}(b).

In the next lemma we recall a
known canonical form of pairs
for feedback similarity. In
the case of controllable
systems, it is known as the
Brunovsky canonical form
\cite{bru}. It can be deduced
from the canonical form of
matrix pencils
\cite[Proposition 3.3]{hin1}.
A much more general canonical
matrix problem was solved in
\cite[\S\,2]{naz_bon}.

Denote by $0_{mn}$ the
$m\times n$ zero matrix;
$0_m:=0_{mm}$. It is agreed
that there exists exactly one
matrix of size $0\times n$
and there exists exactly one
matrix of size $n\times 0$
for every nonnegative integer
$n$; they give the linear
mappings ${\mathbb C}^n\to 0$
and $0\to {\mathbb C}^n$ and
are considered as zero
matrices $0_{0n}$ and
$0_{n0}$. For any $p$-by-$q$
matrix $M_{pq}$, we have
$$
M_{pq}\oplus 0_{m0}
=\begin{bmatrix}
   M_{pq} & 0 \\
  0 &0_{m0}
\end{bmatrix}=\begin{bmatrix}
   M_{pq} & 0_{p0} \\
   0_{mq} & 0_{m0}
\end{bmatrix}=\begin{bmatrix}
 M_{pq}\\0_{mq}
\end{bmatrix}
$$
and
$$
M_{pq}\oplus
0_{0n}=\begin{bmatrix}
  M_{pq} & 0 \\
  0 & 0_{0n}
\end{bmatrix}=\begin{bmatrix}
  M_{pq} & 0_{pn} \\
  0_{0q} & 0_{0n}
\end{bmatrix}=\begin{bmatrix}
   M_{pq} & 0_{pn}
\end{bmatrix}.
$$
Denote
$$
J_k(\lambda):=\begin{bmatrix}
 \lambda&&&0\\1&\lambda&&\\&\ddots &\ddots
 &\\0&&1&\lambda
 \end{bmatrix}\qquad\text{($k$-by-$k$)},
$$
\begin{equation}\label{2.0}
  F_{rl}:=\begin{bmatrix}
 I_r&0_{r,l-r}
 \end{bmatrix},\qquad G_{rl}:=
\begin{bmatrix}
 0_{r,l-r}&I_r
 \end{bmatrix}\quad (0\le r\le l);
\end{equation}
in particular,
$F_{0n}=G_{0n}=0_{0n}$. The
{\it direct sum} of matrix
$t$-tuples is defined as
follows:
\[
(A_1,\dots,A_t)\oplus
(B_1,\dots,B_t):=(A_1\oplus
B_1 ,\dots,A_t\oplus B_t).
\]

\begin{lemma} \label{t3.1}
Every $m\times (n,m)$ pair
$(B,A)$ is feedback similar
to a direct sum of pairs of
the form
\begin{equation}\label{2.1a}
  \left([1\, 0\,\dots\, 0]^T,
  J_k(0)\right),\qquad
 (0_{k0},J_k(\lambda)),\qquad
  (0_{01},0_{0}).
\end{equation}
This sum is determined by
$(B,A)$ uniquely up to
permutation of summands.
\end{lemma}

\begin{proof}
Let $(B,A)$ be an $m\times
(n,m)$ pair. If $B=0$, then
\[
(B,A)=(0_{m0},A)\oplus
(0_{0n},0_{00})=
(0_{m0},A)\oplus
(0_{01},0_{00})\oplus
\dots\oplus (0_{01},0_{00}).
\]
The summand $(0_{m0},A)$ is
feedback similar to a direct
sum of pairs of the form
$(0_{k0},J_k(\lambda))$.

Suppose $B\ne 0$. Then
$(B,A)$ reduces by
transformations (i$'$) and
(ii$'$) to the form
\[
\left(\begin{bmatrix}
 I_r&0\\0&0
 \end{bmatrix},\
 C\right)=\left(\begin{bmatrix}
 I_r\\0
 \end{bmatrix},\
 C\right)\oplus (0_{01},0_{0})
 \oplus
\dots\oplus (0_{01},0_{0}).
\]
The first summand reduces by
transformations (iii$'$) to
the form
\begin{equation}\label{3.a}
{\cal H}(M,N):=
\left(\begin{bmatrix}
 I_r\\0
 \end{bmatrix},\
 \begin{bmatrix}
 0_{r}&0_{r,m-r}\\M&N
 \end{bmatrix}\right).
\end{equation}

If $(M,N)$ is feedback
similar to $(M',N')$, that
is,
\begin{equation}\label{jyk}
S\,[M'\; N']=
  [M\; N]\begin{bmatrix}
    P&U\\0&S
  \end{bmatrix},
\end{equation}
then
\begin{equation}\label{lki}
\begin{bmatrix}
    P&U\\0&S
  \end{bmatrix}
\begin{bmatrix}
 I_k&0&0\\0&M'&N'
\end{bmatrix}
 =\begin{bmatrix}
 I_k&0&0\\0&M&N
\end{bmatrix}
  \begin{bmatrix}
    P&UM'&UN'\\
    0&P&U\\0&0&S
  \end{bmatrix},
\end{equation}
and so ${\cal H}(M,N)$ is
feedback similar to ${\cal
H}(M',N')$. Using induction
in $m+n$, we may assume that
$(M,N)$ is feedback similar
to a direct sum of pairs
${\cal P}_i$ of the form
\eqref{2.1a}. Then ${\cal
H}(M,N)$ is feedback similar
to the direct sum of the
pairs
\[
{\cal H}({\cal P}_i)=
  \begin{cases}
\left([1\, 0\,\dots\,
0]^T,J_{k+1}(0)\right)
    & \text{if ${\cal P}_i=
    \left([1\, 0\,\dots\, 0]^T,
  J_k(0)\right)$}, \\
{\cal P}_i
    & \text{if ${\cal P}_i
    =(0_{k0},J_k(\lambda))$}, \\
(I_1,0_{1}) & \text{if ${\cal
P}_i=(0_{01},0_{0})$.}
  \end{cases}
\]
The uniqueness of this
decomposition follows, for
instance, from \cite[Theorem
2.2]{ser}, in which the
uniqueness of decompositionss
into indecomposables is
proved for all linear matrix
problems.
\end{proof}

For each $m\times (n,m)$ pair
${\cal P}=(B,A)$, we define
the block matrix
\begin{equation}\label{krd}
[{\cal P}]:=[B\ A].
\end{equation}
Consider a category, whose
objects are $m\times (n,m)$
pairs and each morphism from
$(B,A)$ to $(B',A')$ is a
matrix triple $(P,U,S)$ such
that
\begin{equation*}\label{1.4amp}
   [B\ A]
  \begin{bmatrix}
    P&U\\0&S
  \end{bmatrix}=
 S[B'\ A'].
\end{equation*}
By \eqref{1.4a}, two pairs
are isomorphic in this
category if and only if they
are feedback similar. The
next theorem gives canonical
pairs for arbitrarily small
perturbations and feedback
similarity and calculates
their endomorphism rings in
this category.

\begin{theorem} \label{t3.2}
{\rm(a)} \emph{(Canonical
pairs)} In the space
${\mathbb C}^{m\times (n,m)}$
of $m\times (n,m)$ matrix
pairs, $n\ge 1$, all pairs
that are feedback similar to
\begin{equation}\label{3.2}
  {\cal F}_{mn}:=\begin{cases}
    (F_{mn},0_{m}) & \text{if $m\le n$},
     \\[2mm]
    \left(
\begin{bmatrix}
 I_n\\0_{m-n,n}
 \end{bmatrix},\
 \begin{bmatrix}
 0_{nm}\\F_{m-n,m}
 \end{bmatrix}
 \right) & \text{if $m>n$}
  \end{cases}
\end{equation}
$(F_{m-n,m}$ is defined in
\eqref{2.0}$)$ form an open
and everywhere dense set,
which is also the set of all
$m\times (n,m)$ pairs that
are rigid with respect to
feedback similarity.

Alternatively, instead of
${\cal F}_{mn}$ one can take
\begin{equation}\label{3.2a}
  {\cal H}_{mn}:=\begin{cases}
    (G_{mn},0_{m}) & \text{if $m\le n$},
     \\[2mm]
    \left(
\begin{bmatrix}
 I_n\\0_{m-n,n}
 \end{bmatrix},\
 \begin{bmatrix}
 0_{nm}\\H_{m-n,m}
 \end{bmatrix}
 \right) & \text{if $m>n$},
  \end{cases}
\end{equation}
where
\begin{equation}\label{htrd}
 H_{m-n,m}:=\begin{bmatrix}
   I_{(\alpha-1)n}&0&0\\
   0&G_{\beta n}&0_{\beta}
 \end{bmatrix},
\end{equation}
$G_{\beta n}$ is defined in
\eqref{2.0}, and $\alpha$ and
$\beta$ are nonnegative
integers defined as follows:
\begin{equation}\label{hui}
m=\alpha n+\beta,\qquad
0<\beta\le n.
\end{equation}

{\rm(b)} \emph{(Endomorphisms
of canonical pairs)} The
equality
\begin{equation}\label{1.4am}
   [{\cal H}_{mn}]
  \begin{bmatrix}
    P&U\\0&S
  \end{bmatrix}=
 S[{\cal H}_{mn}],\qquad
 P\in{\mathbb C}^{n\times
 n},\ S\in{\mathbb C}^{m\times
 m},\ U\in{\mathbb C}^{n\times
 m},
\end{equation}
$($see \eqref{1.4a} and
\eqref{krd}$)$ holds if and
only if for some $S_1\in
{\mathbb C}^{(n-\beta)\times
(n-\beta)}$, $S_3\in {\mathbb
C}^{\beta\times\beta}$, and
$S_2,S_4\in {\mathbb
C}^{(n-\beta)\times\beta}$ we
have
\begin{equation}\label{ldf}
\begin{bmatrix}
    P&U\\0&S
  \end{bmatrix}=
  R_{\alpha+2}(S_1,S_2,S_3,S_4),
\qquad
S=R_{\alpha+1}(S_1,S_2,S_3,S_4),
\end{equation}
where
\begin{align}\nonumber
&\
\overbrace{\quad}^{n-\beta} \
\overbrace{\quad}^{\beta}\
   \qquad\qquad\qquad\qquad\ \,
 \overbrace{\quad}^{n-\beta}
\:
\overbrace{\quad}^{\beta}\:
\overbrace{\quad}^{\beta}
            \\[-2mm] \label{4.5a}
 R_{\gamma}(S_1,S_2,S_3,S_4):=&
 \left[
\begin{tabular}{cc|cc|cc|cc|c}
 $S_1$&$S_2$   & $0$&$S_4$ &&&&\\
 &$S_3$ & $0$&$0$ &&&& \\ \hline
 && $\ddots$& & $\ddots$& & &&
 \\[-1.5mm]
 && &$\ddots$ & &$\ddots$ & && \\ \hline
 &&&&   $S_1$&$S_2$   & $0$&$S_4$ \\
 &&&& &$S_3$ & $0$&$0$  \\ \hline
 &&&&&& $S_1$&$S_2$ & $S_4$ \\
 &&&&&& &$S_3$ & $0$  \\ \hline
 &&&&&& & & $S_3$
\end{tabular}\right]
\end{align}
$(\gamma$ is the number of
diagonal blocks $S_3$;
unspecified blocks are
zero$)$.

In this statement one can
replace ${\cal H}_{mn}$ by
${\cal F}_{mn}$, which is
simpler, but then
$R_{\gamma}(S_1,\dots,S_4)$
must be replaced by
\begin{equation}\label{krf}
 \left[
\begin{tabular}{cc|cc|cc|c}
 $S_1$&$0$   & $0$&$0$ &&\\
 $S_2$&$S_3$ & $S_4$&$0$ && \\ \hline
 && $\ddots$& & $\ddots$& &
 \\[-1.5mm]
 && &$\ddots$ & &$\ddots$ &  \\ \hline
 &&&& $S_1$&$0$ & $0$ \\
 &&&& $S_2$&$S_3$ & $S_4$  \\ \hline
 &&&& & & $S_3$
\end{tabular}\right],
\end{equation}
which is not
block-triangular.
\end{theorem}

\begin{proof}
(a) Let $(B,A)\in {\mathbb
C}^{m\times (n,m)}$, $n\ge
1$. First, we make $\rank B=
\min(m,n)$ by an arbitrarily
small perturbation and reduce
$B$ to the form
\begin{equation*}\label{2.6'}
  \begin{cases}
   F_{mn}=\begin{bmatrix}
  I_m&0
\end{bmatrix} & \text{if
$m\le n$}, \\[2mm]
  \begin{bmatrix}
  I_n\\0
\end{bmatrix}
 & \text{if $m> n$,}
  \end{cases}
\end{equation*}
using transformations (i$'$)
and (ii$'$). Then we reduce
the pair by transformations
(iii$'$) to the form
$(F_{mn},0_{m})$ if $m\le n$
or to the form \eqref{3.a}
with $r=n$ if $m>n$.

If $m>n$, using induction in
$m$ we can assume that
$(M,N)$ reduces by an
arbitrary small perturbation
to some $(M+\Delta M,
N+\Delta N)$ being feedback
similar to ${\cal F}_{m-n,n}$
defined in \eqref{3.2}. By
\eqref{jyk} and \eqref{lki},
${\cal H}(M+\Delta M,
N+\Delta N)$ is feedback
similar to ${\cal F}_{mn}$.
Reasoning as in Corollary
\ref{corrd}, we can prove
that $(B,A)$ is transformed
by an arbitrarily small
perturbation to a pair that
is feedback similar to ${\cal
F}_{mn}$. Hence, the set
$\cal S$ of pairs that are
feedback similar to ${\cal
F}_{mn}$ is everywhere dense
in ${\mathbb C}^{m\times
(n,m)}$. Since ${\cal
F}_{mn}$ is rigid, there
exists its neighborhood $V$
in ${\mathbb C}^{m\times
(n,m)}$ such that
$V\subset\cal S$. For any
pair ${\cal P}\in{\cal S}$,
there is a transformation of
feedback similarity that
transforms ${\cal F}_{mn}$ to
$\cal P$; it also transforms
$V$ to some neighborhood $W$
of $\cal P$. Since each pair
in $V$ is feedback similar to
${\cal F}_{mn}$, each pair in
$W$ is also feedback similar
to ${\cal F}_{mn}$, hence
$W\subset {\cal S}$.
Therefore, each pair ${\cal
P}\in{\cal S}$ possesses a
neighborhood that is
contained in ${\cal S}$, and
so the set ${\cal S}$ is open
in ${\mathbb C}^{m\times
(n,m)}$.

If $m>n$, then the pair
${\cal F}_{mn}$ reduces to
${\cal H}_{mn}$ in
\eqref{3.2a} by those
permutations of rows and
columns that are special
cases of transformations
(i$'$) and (ii$'$).

(b) Assume first that
\eqref{1.4am} holds. If $m\le
n$, then $[{\cal
H}_{mn}]=[G_{mn}\
0_m]=[[0_{m,n-m}\ I_m]\
0_m]$, and so \eqref{1.4am}
ensures
\[
\begin{bmatrix}
    P&U\\0&S
  \end{bmatrix}=
\begin{bmatrix}
    \begin{bmatrix}
    *&*\\0&S
  \end{bmatrix}&\begin{bmatrix}
    *\\0
  \end{bmatrix}\\0&S
  \end{bmatrix}=R_2(*,*,S,*).
\]
We have \eqref{ldf} since
$R_1(*,*,S,*)=S$ and by
\eqref{hui} $\alpha=0$.

Let $m>n$. Equating the
corresponding vertical strips
in \eqref{1.4am}, we obtain
\begin{equation}\label{4.3f}
\begin{bmatrix}
 P\\0_{m-n,n}
 \end{bmatrix}
=S\begin{bmatrix}
  I_n\\0_{m-n,n}
 \end{bmatrix},
                    \qquad
 \begin{bmatrix}
  U\\H_{m-n,m}S
 \end{bmatrix}
=
 S  \begin{bmatrix}
 0_{nm}\\H_{m-n,m}
 \end{bmatrix}.
\end{equation}

Let us prove that
$S=R_{\alpha+1}(S_1,S_2,S_3,S_4)$
for some $S_1,\dots,S_4$.
Partition $S$ into blocks
\[
S=\begin{bmatrix}
 S_{11}&\dots &S_{1,\alpha+1}\\
 \hdotsfor{3}\\
 S_{\alpha+1,1}&\dots& S_{\alpha+1,\alpha+1}
 \end{bmatrix}
\]
with $n\times n,\dots,n\times
n, \beta\times \beta$
diagonal blocks. By the first
equality in \eqref{4.3f},
\begin{equation}\label{4.5f}
 P=S_{11},\qquad
 S_{21}=\dots =S_{\alpha 1}=0,\qquad
  S_{\alpha+1,1}=0.
\end{equation}
Since
\[
H_{m-n,m}=\begin{bmatrix}
 I_n&&&0&0_{n\beta}
 \\&\ddots&&&\vdots
 \\&&I_n&&0_{n\beta}
 \\0&&&G_{\beta n}&0_{\beta}
 \end{bmatrix},\qquad G_{\beta n}=
\begin{bmatrix}
 0_{\beta,n-\beta}&I_{\beta}
 \end{bmatrix},
\]
by the second equality in
\eqref{4.3f} we have
\begin{equation}\label{lev}
U=\begin{bmatrix}
    S_{12}&\dots&S_{1\alpha}
    &S_{1,\alpha+1}
    G_{\beta n}&0
  \end{bmatrix}
\end{equation}
and
\begin{multline*}\label{4.5aa}
  \begin{bmatrix}
    S_{11}&\dots&S_{1,\alpha-1}&S_{1\alpha}
    &S_{1,\alpha+1}\\
    \vdots&\ddots&\vdots&\vdots&\vdots\\
    S_{\alpha-1,1}&\dots&
    S_{\alpha-1,\alpha-1}&S_{\alpha-1,\alpha}
    &S_{\alpha-1,\alpha+1}\\
    G_{\beta n}S_{\alpha 1}&\dots&
    G_{\beta n}S_{\alpha,\alpha-1}&
    G_{\beta n}S_{\alpha\alpha}
    &G_{\beta n}S_{\alpha,\alpha+1}
  \end{bmatrix}
 \\=
  \begin{bmatrix}
    S_{22}&\dots&S_{2\alpha}
    &S_{2,\alpha+1}
    G_{\beta n}&0\\
 \vdots&\ddots&\vdots&\vdots&\vdots\\
    S_{\alpha 2}&\dots&
    S_{\alpha\alpha}
   &S_{\alpha,\alpha+1}
   G_{\beta n}&0\\
    S_{\alpha+1,2}&\dots&
    S_{\alpha+1,\alpha}
    &S_{\alpha+1,\alpha+1}
    G_{\beta n}&0
  \end{bmatrix}.
\end{multline*}
Let us equate the entries of
these matrices along each
line that is parallel with
the main diagonal:
\begin{itemize}
  \item[(a)]
The equalities
\begin{align*}
 S_{1,\alpha+1}&=0,\\
 S_{1\alpha}=S_{2,\alpha+1}G_{\beta n},\qquad
 S_{2,\alpha+1}&=0,\\
 S_{1,\alpha-1}=S_{2,\alpha}=
 S_{3,\alpha+1}G_{\beta n},\qquad
 S_{3,\alpha+1}&=0,\\
      &\vdots\\
 S_{13}=S_{24}=\dots=
 S_{\alpha-2,\alpha}=
 S_{\alpha-1,\alpha+1}G_{\beta n},\qquad
 S_{\alpha-1,\alpha+1}&=0
\end{align*}
imply $S_{ij}=0$ if $j-i\ge
2$.

  \item[(b)]
The equalities
$$
 S_{12}=S_{23}=\dots=S_{\alpha-1,\alpha}=
 S_{\alpha,\alpha+1}G_{\beta n},\qquad
 G_{\beta n} S_{\alpha,\alpha+1}=0
$$
imply
$$
 S_{12}=S_{23}=\dots=S_{\alpha-1,\alpha}=
\begin{bmatrix}
  0&S_4\\0&0
\end{bmatrix},\qquad
 S_{\alpha,\alpha+1}=\begin{bmatrix}
  S_4\\0
\end{bmatrix}
$$
for some
$(n-\beta)\times\beta$ matrix
$S_4$.

  \item[(c)]
The equalities
$$
 S_{11}=S_{22}=\dots=S_{\alpha\alpha},
 \qquad
 G_{\beta n}S_{\alpha\alpha}=
 S_{\alpha+1,\alpha+1}G_{\beta n}
$$
imply
\begin{equation*}\label{4.6aaf}
  S_{11}=S_{22}=\dots=S_{\alpha\alpha}=
\begin{bmatrix}
  S_1&S_2\\0&S_3
\end{bmatrix},\qquad S_3:=
S_{\alpha+1,\alpha+1}.
\end{equation*}

  \item[(d)]
The equalities
\begin{align*}
 S_{21}=S_{32}=\dots=S_{\alpha,\alpha-1},
 \qquad G_{\beta n}S_{\alpha,\alpha-1}
 &=S_{\alpha+1,\alpha},\\
 S_{31}=\dots=S_{\alpha,\alpha-2},
 \qquad G_{\beta n}S_{\alpha,\alpha-2}
 &=S_{\alpha+1,\alpha-1},\\
   &\vdots\\
 S_{\alpha-1,1}=S_{\alpha 2},
 \qquad G_{\beta n}S_{\alpha 2}
 &=S_{\alpha+1,3},\\
 \qquad G_{\beta n}S_{\alpha 1}
 &=S_{\alpha+1,2}
\end{align*}
and \eqref{4.5f} imply
$S_{ij}=0$ if $i>j$.
\end{itemize}
This proves the first
equality in \eqref{ldf}. The
second equality in
\eqref{ldf} follows from
\eqref{lev} and the first
equality in \eqref{4.5f}.

Conversely, the equalities
\eqref{ldf} ensure
\eqref{1.4am}. For example,
if $\alpha=1$, then
\eqref{1.4am} takes the form
\begin{multline}\label{kex}
\begin{bmatrix}
I_n&0_{n}&0
 \\0&[0\ I_{\beta}]&0_{\beta}
 \end{bmatrix}
\begin{bmatrix}
\begin{bmatrix}
  S_1&S_2\\0&S_3
\end{bmatrix}&\begin{bmatrix}
  0&S_4\\0&0
\end{bmatrix}&0
 \\0&\begin{bmatrix}
  S_1&S_2\\0&S_3
\end{bmatrix}&\begin{bmatrix}
  S_4\\0
\end{bmatrix}\\0&0&S_3
 \end{bmatrix}
   \\
=\begin{bmatrix}
\begin{bmatrix}
  S_1&S_2\\0&S_3
\end{bmatrix}&\begin{bmatrix}
  S_4\\0
\end{bmatrix}\\0&S_3
 \end{bmatrix}
\begin{bmatrix}
I_n&0_{n}&0
 \\0&[0\ I_{\beta}]&0_{\beta}
 \end{bmatrix}
\end{multline}
\end{proof}

\begin{remark}\label{rem}
The condition $n\ge 1$ in
Theorem \ref{t3.2}(a) is
essential: each
$m\times(0,m)$ pair is
transformed by an arbitrarily
small perturbation to a pair
that is feedback similar to
$(0_{m0},\diag(\lambda_1,\dots,\lambda_m))$
with distinct
$\lambda_1,\dots,\lambda_m$
determined up to small
perturbations.
\end{remark}

\section{Feedback similarity of triples}
\label{s4}

The next lemma is proved by
using several steps of
Belitskii's algorithm
\cite{bel,ser} and
arbitrarily small
perturbations.

\begin{lemma} \label{t4.1}
Every $m\times (n,m,m)$
triple $(C,B,A)$, $n\ge 1$,
is transformed by an
arbitrarily small
perturbation to a triple that
is feedback similar to
\begin{equation}\label{4.1}
  \begin{cases}
(G_{mn},0,0) & \text{if $m\le n$}, \\
 {\cal K}(N):=\left(
\begin{bmatrix}
 I_n\\0_{m-n,n}
 \end{bmatrix},
 \begin{bmatrix}
 0_{nm}\\ H_{m-n,m}
 \end{bmatrix},
 \begin{bmatrix}
 0_{nm}\\N
 \end{bmatrix}\right)
 & \text{if $m>n$},
  \end{cases}
\end{equation}
where $N$ is some
$(m-n)\times m$ matrix and $
H_{m-n,m}$ is defined in
\eqref{htrd}.

Two triples ${\cal K}(N)$ and
${\cal K}(N')$ are feedback
similar if and only if
\begin{equation}\label{4.2}
 N'= R_{\alpha}(S_1,S_2,S_3,S_4)^{-1}\cdot
 N\cdot
 R_{\alpha+1}(S_1,S_2,S_3,S_4)
\end{equation}
$($see \eqref{ldf}$)$ for
some $S_2,S_4\in{\mathbb
C}^{(n-\beta)\times\beta}$
and nonsingular matrices
$S_1\in{\mathbb
C}^{(n-\beta)\times
(n-\beta)}$ and
$S_3\in{\mathbb
C}^{\beta\times\beta}$.
\end{lemma}

\begin{proof}
Let $(C,B,A)$ be an $m\times
(n,m,m)$ triple, $n\ge 1$. By
Theorem \ref{t3.2}(a), there
is an arbitrarily small
perturbation of $(C,B)$ such
that the obtained pair
$(C+\Delta C, B+\Delta B)$ is
feedback similar to the pair
${\cal H}_{mn}$ in
\eqref{3.2a}, and then
$(C+\Delta C, B+\Delta B,A)$
is feedback similar to
\eqref{4.1}.

Let $N,N'\in{\mathbb
C}^{(m-n)\times m}$. Suppose
first that ${\cal K}(N)$ and
${\cal K}(N')$ are feedback
similar. By \eqref{1.4b},
\begin{equation}\label{4.3a}
\begin{bmatrix}
 I_n&0&0\\0&H_{m-n,m}&N
 \end{bmatrix}
 \begin{bmatrix}
    P&U&V\\0&S&0\\0&0&S
  \end{bmatrix}
=
 S \begin{bmatrix}
 I_n&0&0\\0&H_{m-n,m}&N'
 \end{bmatrix}
\end{equation}
for some $U,V$ and
nonsingular $P$ and $S$. Then
\eqref{1.4am} holds, which
ensures \eqref{ldf}. Equating
the last vertical strips of
the matrices in \eqref{4.3a}
gives
\begin{equation*}\label{4.3}
 \begin{bmatrix}
  V\\ NS
 \end{bmatrix}=
 S\begin{bmatrix}
    0_{nm}\\N'
  \end{bmatrix},
\end{equation*}
which defines $V$ and ensures
\eqref{4.2}.

Conversely, if \eqref{4.2}
holds, then by analogy with
\eqref{kex} we have
\eqref{4.3a} for
\[
P=\begin{bmatrix}
  S_1&S_2\\0&S_3
\end{bmatrix},\qquad
U=\begin{bmatrix}
  0&S_4&0&\dots&0\\0&0&0&\dots&0
\end{bmatrix},\qquad
V=UN',
\]
and
$S=R_{\alpha+1}(S_1,S_2,S_3,S_4)$.
Hence, ${\cal K}(N)$ and
${\cal K}(N')$ are feedback
similar.
\end{proof}

\begin{remark}
Instead of $G_{mn}$ and
$H_{m-n,m}$ in \eqref{4.1},
one may take $F_{mn}$ and
$F_{m-n,m}$ replacing in
\eqref{4.2} the matrix $
R_{\gamma}(S_1,\dots,S_4)$
defined in \eqref{4.5a} with
\eqref{krf}. We prefer
\eqref{4.1} since the matrix
\eqref{4.5a} is upper
block-triangular and we can
reduce $N$ to Belitskii's
canonical form \cite{bel,ser}
by transformations
\eqref{4.2} preserving the
other blocks of ${\cal
K}(N)$. Examples of this
reduction are given in
Theorems \ref{c4.2} and
\ref{c4.1}.
\end{remark}

\begin{theorem}\label{c4.2}
Each $m\times (1,m,m)$ triple
$(C,B,A)$, $m\ge 2$, reduces
by an arbitrarily small
perturbation to a triple that
is feedback similar to
\[
  \left(
\begin{bmatrix}
 1\\0\\ \vdots \\0
 \end{bmatrix},\
 \begin{bmatrix}
 0&\dots&0&0\\1&\dots&0&0\\
 \vdots&\ddots&\vdots&\vdots\\
 0&\dots&1&0
 \end{bmatrix},\
 \begin{bmatrix}
 0&\dots&0&0\\ *&\dots
 &*&*\\
 \vdots&\ddots&\vdots&\vdots\\
 *&\dots
 &*&*
 \end{bmatrix}\right)
\]
$($the stars denote
unspecified entries$)$. This
triple is determined by
$(C,B,A)$ uniquely up to
small perturbations of the
entries denoted by stars.

In greater generality, each
$\alpha n\times (n,\alpha
n,\alpha n)$ triple
$(C,B,A)$, $\alpha\ge 2$,
reduces by an arbitrarily
small perturbation to a
triple that is feedback
similar to a triple of the
form
\begin{equation}\label{4.7}
  \left(
\begin{bmatrix}
 I_n\\0\\ \vdots \\0
 \end{bmatrix},\
 \begin{bmatrix}
 0&\dots&0&0\\I_n&\dots&0&0\\
 \vdots&\ddots&\vdots&\vdots\\
 0&\dots&I_n&0
 \end{bmatrix},\
 \begin{bmatrix}
 0&\dots&0&0\\N_{11}&\dots
 &N_{1,\alpha-1}&N_{1\alpha}\\
 \vdots&\ddots&\vdots&\vdots\\
 N_{\alpha-1, 1}&\dots
 &N_{\alpha-1,\alpha-1}&N_{\alpha-1,\alpha}
 \end{bmatrix}\right),
\end{equation}
in which all blocks are
$n$-by-$n$,
\[
  N_{11}=\diag(
    \lambda_1,\lambda_2,\dots,
    \lambda_n)\qquad
(\text{$\lambda_1,\dots,\lambda_n$
are distinct}),
\]
\[
 N_{12}=\begin{bmatrix}
    *&1&\dots&1 \\
    *&*&\dots&*\\
    \vdots&\vdots&\ddots&\vdots\\
    *&*&\dots&*
  \end{bmatrix}\qquad
(\text{the stars denote
unspecified entries}),
\]
and the other $N_{ij}$ are
arbitrary. The triple
\eqref{4.7} is determined by
$(C,B,A)$ uniquely up to
small perturbations of
$\lambda_1,\dots,\lambda_n$
in $N_{11}$, of the entries
denoted by stars in $N_{12}$,
and of the entries in the
other $N_{ij}$.
\end{theorem}

\begin{proof}
Let $(C,B,A)$ be $\alpha
n\times (n,\alpha n,\alpha
n)$, $\alpha\ge 2$. By Lemma
\ref{t4.1}, $(C,B,A)$ reduces
by an arbitrarily small
perturbation and a feedback
similarity transformation to
a triple of the form
\eqref{4.7}, in which
$N_{ij}$ are $n$-by-$n$. We
can reduce $N:=[N_{ij}]$ by
transformations \eqref{4.2}
preserving the other blocks
of the triple \eqref{4.7}. We
have
\begin{equation*}\label{4.8}
 R_{\gamma}(S_1,S_2,S_3,S_4)
 =S_3\oplus\dots\oplus S_3\qquad
 \text{(${\gamma}$
 summands)}
\end{equation*}
since $\beta=n$ in
\eqref{hui} and so $S_1$ is
$0\times 0$ in \eqref{4.5a}.
Hence we can reduce all
$N_{ij}$ by simultaneous
similarity transformations
\begin{equation}\label{4.9}
 N_{ij}'=S_3N_{ij}S_3^{-1},\qquad
 1\le i\le \alpha-1,\quad
 1\le j\le \alpha.
\end{equation}
By an arbitrarily small
perturbation and some
transformation \eqref{4.9} we
reduce $N_{11}$ to a diagonal
matrix with distinct diagonal
entries. To preserve $N_{11}$
we must reduce the other
blocks $N_{ij}$ by
transformations \eqref{4.9}
with diagonal $S_3$. Using an
arbitrarily small
perturbation we make nonzero
the $(1,2),\dots,(1,n)$
entries of the first row of
$N_{12}$ and reduce them to 1
by transformations
\eqref{4.9} with diagonal
$S_3$. Each transformation
\eqref{4.9} that preserves
$N_{11}$ and the
$(1,2),\dots,(1,n)$ entries
of $N_{12}$ is the identity
transformation, so we can
reduce the other entries of
$N_{ij}$ only by arbitrarily
small perturbations.
\end{proof}

For every $p\times (q,p,p)$
triple $(C,B,A)$, we define
the $(2p+q)\times
(p+q,2p+q,2p+q)$ triple
\begin{equation*}\label{4.9a}
{\cal L}(C,B,A):= \left(
\begin{bmatrix}
 I_{p+q}\\0_{p,p+q}
 \end{bmatrix},\
 \begin{bmatrix}
 0&0&0\\
 0_{pq}&I_p&0_{p}
  \end{bmatrix},\
 \begin{bmatrix}
 0&0&0\\
 C&B&A
 \end{bmatrix}\right).
\end{equation*}
Put
\[
{\cal L}^{(i)}(C,B,A):=
\underbrace{{\cal
L}\dots{\cal
L}}_{\mbox{$i$-times}}(C,B,A),\qquad
i=0,1,2,\ldots.
\]

\begin{theorem}\label{c4.1}
Let $m$ and $n$ be natural
numbers, and let ${\mathbb
C}^{m\times (n,m,m)}$ denote
the space of all $m\times
(n,m,m)$ triples.

{\rm(a)} If
$m<n(1+\sqrt{5})/{2}$, then
there is exactly one $m\times
(n,m,m)$ triple of the form
${\cal
L}^{(l)}(F_{pq},0_{p},0_{p})$.
All triples that are feedback
similar to it form an open
and everywhere dense set in
${\mathbb C}^{m\times
(n,m,m)}$, which coincides
with the set of all $m\times
(n,m,m)$ triples that are
rigid with respect to
feedback similarity.

{\rm(b)} If
$m>n(1+\sqrt{5})/{2}$, then
all $m\times (n,m,m)$ triples
are not rigid with respect to
feedback similarity.
\end{theorem}

\begin{proof}
Let $m$ and $n$ be natural
numbers, and let $(C,B,A)$ be
${m\times (n,m,m)}$. We say
that a triple ${\cal T}$
\emph{reduces} to a triple
${\cal T}^{\,\prime}$ if
${\cal T}$ reduces to ${\cal
T}^{\,\prime}$ by an
arbitrarily small
perturbation and a feedback
similarity transformation.
\medskip

(a) Suppose first that
\begin{equation}\label{ked}
m<n(1+\sqrt{5})/{2}\approx
1.618 n,
\end{equation}
and prove by induction on
$m-n$ that $(C,B,A)$ reduces
to some ${\cal
L}^{(l)}(F_{pq},0_{p},0_{p})$.

The base of induction is
trivial: if $m\le n$, then by
Lemma \ref{t4.1} $(C,B,A)$
reduces to
$(F_{mn},0_{m},0_{m})$. This
triple is rigid for feedback
similarity and is feedback
similar to each rigid
$m\times (n,m,m)$ triple.

Let $m>n$. Then by
\eqref{hui} and \eqref{ked}
we have $\alpha=1$ and
$\beta=m-n$. According to
Lemma \ref{t4.1}, $(C,B,A)$
reduces to some triple
\begin{align*}\label{mep}
 {\cal K}([C'\ B'\ C'])&=
 \left(
\begin{bmatrix}
 I_n\\0_{\beta n}
 \end{bmatrix},
 \begin{bmatrix}
0_{n,n-\beta}&0_{n\beta}&0_{n\beta}
\\ 0_{\beta,n-\beta}
&I_{\beta}&0_{\beta}
 \end{bmatrix},
 \begin{bmatrix}
0_{n,n-\beta}&0_{n\beta}&0_{n\beta}
\\ C'&B'&C'
 \end{bmatrix}\right)\\
 &={\cal L}(C',B',A'),
\end{align*}
in which $(C',B',A')$ is
$m'\times (n'\times m'\times
m')$ and
\begin{equation}\label{4.11}
 m':=m-n,\qquad n':=-m+2n.
\end{equation}
By Lemma \ref{t4.1}, ${\cal
K}([C'\ B'\ C'])$ is feedback
similar to ${\cal K}([C'_1\
B'_1\ C'_1])$ if and only if
there exist $U,V$, and
nonsingular $P$ and $S$ such
that
\begin{equation}\label{4.12}
  [C'_1\;B'_1\;A'_1]=
  S^{-1}[C'\;B'\;A']
  \begin{bmatrix}
    P&U&V\\0&S&0\\0&0&S
   \end{bmatrix}
 \end{equation}
(the last matrix is
$R_2(P,U,S,V)$ defined in
\eqref{4.5a}). Therefore,
${\cal L}(C',B',A')$  is
feedback similar to ${\cal
L}(C'_1,B'_1,A'_1)$ if and
only if $(C',B',A')$ is
feedback similar to
$(C'_1,B'_1,A'_1)$.

The numbers $m'$ and $n'$ are
natural: $m'>0$ since $m>n$,
and $n'=2n-m>0$ since
$1.7n-m>0$ by \eqref{ked}.
Furthermore, $m'<
n'(1+\sqrt{5})/2$ because by
\eqref{ked}
$$
\frac{m_1}{n_1}=\frac{m-n}{-m+2n}
 <\frac{n(1+\sqrt{5})/2-n}{
 -n(1+\sqrt{5})/2+2n}
 =\frac{-1+\sqrt{5}}{3-\sqrt{5}}
 =\frac{1+\sqrt{5}}{2}.
$$
Since $(m-n)-(m'-n')=n'>0$,
the induction hypothesis
ensures that $(C',B',A')$
reduces to ${\cal
L}^{(l-1)}(F_{pq},0_{p},0_{p})$
for some $p,q$, and $l$ that
are uniquely determined by
$m'$ and $n'$. Therefore,
$(C,B,A)$ reduces to ${\cal
K}([C'\ B'\ C']) ={\cal
L}(C',B',A')$, which reduces
to ${\cal
L}^{(l)}(F_{pq},0_{p},0_{p})$
that is uniquely determined
by $m$ and $n$.

We have proved that all
${m\times (n,m,m)}$ triples
reduce to the same triple
${\cal L}:={\cal
L}^{(l)}(F_{pq},0_{p},0_{p})$,
and so the set $\cal S$ of
all triples that are feedback
similar to ${\cal L}$ is
everywhere dense. Since
${\cal L}$ is rigid with
respect to feedback
similarity, there exists its
neighborhood $V$ that is
contained in $\cal S$. For
any triple ${\cal T}\in{\cal
S}$, there is a
transformation of feedback
similarity that transforms
${\cal L}$ to $\cal T$; it
also transforms $V$ to some
neighborhood $W$ of $\cal T$.
Since each triple in $V$ is
feedback similar to ${\cal
L}$, each triple in $W$ is
also feedback similar to
${\cal L}$, hence $W\subset
{\cal S}$. Therefore, each
triple ${\cal T}\in{\cal S}$
possesses a neighborhood that
is contained in ${\cal S}$,
and so the set ${\cal S}$ is
open. {}\medskip

(b) Let
\begin{equation}\label{msi}
m\ge n(1+\sqrt{5})/2.
\end{equation}
Since $(C,B,A)$ is fixed, the
equality \eqref{1.4b} defines
the mapping
\[
f\colon {\cal U}\to {\mathbb
C}^{m\times (n,m,m)},\qquad
(S,P,U,V)\longmapsto
(C',B',A'),
\]
where
\[
{\cal
U}:=\{(S,P,U,V)\in{\mathbb
C}^{m\times m} \times
{\mathbb C}^{n\times n}\times
{\mathbb C}^{n\times m}\times
{\mathbb C}^{n\times m}\,|\,
\det(S)\det(P)\ne 0\}.
\]
This mapping is rational
since by \eqref{1.4b} the
entries of $C',\ B'$, and
$A'$ are polynomials (in
entries of $S,P,U$, and $V$)
divided by $\det(S)$. Its
image is the set of all
triples that are feedback
similar to $(C,B,A)$.

Suppose that $(C,B,A)$ is
rigid. Then the image of $f$
contains a neighborhood of
$(C,B,A)$, hence ${\mathbb
C}^{m\times
(n,m,m)}\smallsetminus\im(f)$
can not be dense in ${\mathbb
C}^{m\times (n,m,m)}$, and so
$\dim({\cal U})\ge
\dim({\mathbb C}^{m\times
(n,m,m)})$ by \cite[Section
3, Proposition 1.2]{hir}.
This means that $
m^2+n^2+2mn\ge mn+2m^2,$
\[
(m/n)^2-m/n-1\le 0, \qquad
m/n<(1+\sqrt{5})/2,
\]
which contradicts to
\eqref{msi}. Therefore, there
are no rigid triples of this
size.
\end{proof}

For each $m\times (n,m,m)$
triple ${\cal T}=(C,B,A)$, we
define the ${m\times (n+2m)}$
polynomial matrix
\begin{equation*}\label{4.15}
 {\cal T}(x,y)=
 [C\ \;xI_m+B\ \;yI_m+A].
\end{equation*}
The next lemma is trivial,
but it can be useful.

\begin{lemma}\label{t4.3}
Two matrix triples $\cal T$
and ${\cal T}^{\,\prime}$ are
feedback similar if and only
if the corresponding
polynomial matrices ${\cal
T}(x,y)$ and ${\cal
T}^{\,\prime}(x,y)$ are
\emph{strictly equivalent};
this means that
\begin{equation}\label{4.14}
 S{\cal T}^{\,\prime}(x,y)=
 {\cal T}(x,y)R
\end{equation}
for some nonsingular complex
matrices $S$ and $R$.
\end{lemma}

\begin{proof}
Let ${\cal T}=(C,B,A)$ and
${\cal
T}^{\,\prime}=(C',B',A')$ be
$m\times (n,m,m)$.

If ${\cal T}$ and ${\cal
T}^{\,\prime}$ are feedback
similar, then there exists a
nonsingular matrix
\begin{equation}\label{4.16}
  R=\begin{bmatrix}
    P&U&V\\0&S&0\\0&0&S
  \end{bmatrix}
\end{equation}
such that
\[
S[C'\;B'\;A']=[C\;B\;A]R.
\]
Since
\[
{\cal
T}(x,y)=[C\;B\;A]+x[0\;I\;0]+
y[0\;0\;I]
\]
and
\[
S[0\;I_m\;0]=[0\;I_m\;0]R,\qquad
S[0\;0\;I_m]=[0\;0\;I_m]R,
\]
we have \eqref{4.14}.

Conversely, let \eqref{4.14}
hold. This polynomial
equality breaks into three
scalar equalities:
\begin{gather*}
S[C'\;B'\;A']=[C\;B\;A]R,\\
S[0\;I_m\;0]=[0\;I_m\;0]R,\qquad
S[0\;0\;I_m]=[0\;0\;I_m]R.
\end{gather*}
By the last two equalities,
the matrix $R$ has the form
\eqref{4.16}. So by the first
equality ${\cal T}$ and
${\cal T}^{\,\prime}$ are
feedback similar.
\end{proof}

\begin{remark}
The authors are grateful to the
reviewer for suggestions and the
following commentaries: We study the
orbits of the action \eqref{1.4b} of
the product of two groups on the
space of matrix triples, which can
be identified with $\mathbb
C^{m(n+2m)}$. Namely, from the left
one has the action of GL$(m;\mathbb
C)$, and from the right one has the
action of the $3\times 3$ block
upper triangular subgroup of GL$(n+
2m;\mathbb C)$ with $(2,3)$ block is
equal to zero. It is known that each
orbit under such an action is a
smooth irreducible semi-affine
variety $V$, i.e. its closure is an
affine irreducible variety
$\overline{V}$, and $V =
\overline{V} \setminus W$, where $W$
is a strict subvariety of $\overline
V$. Moreover, all singular points of
$\overline{V}$ are contained in $W$.
The orbits of the maximal dimension
$d$ are called the ``generic''
orbits. Theorem 10 gives the unique
canonical form of a generic orbit.
The parameter space of such orbits
is $m(n + 2m)-d$ dimensional. The
notion of rigid system is equivalent
to the assumption of the existence
of orbits of dimension $m(n + 2m)$.
Since such an orbit $V$ is an
irreducible semi-affine variety, it
follows that $\overline V =\mathbb
C^ {m(n+2m)}$. Hence there is only
one orbit like that as Theorem 11
claims.
\end{remark}

\section*{Acknowledgement}

The authors would like to
thank Professor Diederich
Hinrichsen for very valuable
suggestions on two manuscript
versions of this article.

\end{document}